%% file: cancel_hal.tex
\documentclass[12pt]{amsart}
\usepackage{epsfig,amssymb,epic,eepic,color}
\usepackage[all]{xy}
 \usepackage[T1]{fontenc}

\newcommand{\be}{\begin{enumerate}}
\newcommand{\ee}{\end{enumerate}}
\newcommand{\bi}{\begin{itemize}}
\newcommand{\ei}{\end{itemize}}

\def\R{\mathbb{R}}
\def\M{\mathbb{M}}

\def\al{\alpha}
\def\be{\beta}

\def\De{\Delta}

\def\ep{\varepsilon}

\def\nd{\noindent}
\def\bull{\hfill$\Box$\\}

\begin{document}
\vskip 1cm

\begin{center}

\medskip

{\sc  A proof of Morse's theorem about the cancellation of critical points}

\vskip .5cm

Fran\c cois Laudenbach

\end{center}
\title{}
\address{Laboratoire de Math\'ematiques Jean Leray,  UMR 6629 du CNRS, Facult\'e des Sciences et Techniques,Universit\'e de Nantes, 2, rue de la Houssini\`ere, F-44322 Nantes cedex 3, France.}
\email{francois.laudenbach@univ-nantes.fr}



\vskip .5cm

\begin{itemize}
{\small 
\item[] {\sc Abstract.} In this note, we give a proof of the famous theorem of M. Morse
dealing with the cancellation of a pair of non-degenerate critical points of a smooth function. Our  proof consists
of a reduction to the one-dimensional case where the question becomes easy to answer.  }\\
\end{itemize}


\thispagestyle{empty}
\vskip .5cm
\vskip .5cm
\input english3.2.tex

\end{document}

%% file: english3.2.tex


\vskip .5cm
\section{Introduction}
\medskip
Let us consider an $n$-dimensional closed manifold $M$ equipped with a Morse function
$f:M\to\R$. Marston Morse devoted two papers to the question of canceling a pair of 
critical points. In \cite{morse1} he considered the possibility of canceling  a local minimum 
(index 0) with a critical point of index 1; and, similarly the cancellation
 of a pair of indices
 $(n-1, n)$; in that paper, a {\it polar} function means 
 a function with no supernumerary local extrema. 
 In \cite{morse2} (see also \cite{hueb})
Morse extended  its cancellation criterion to the case of a pair of indices $(k,k+1)$.
It is worth noticing that these two papers show the first time where Morse is considering 
gradients globally in the manifold $M$, in the spirit of the 1949 Note by Ren\'e Thom \cite{thom49}.

S. Smale generalized Morse's criterion and got a more algebraic criterion for cancellation.
This generalization became an important tool in Smale's proof of the {\it generalized Poincar\'e
conjecture} \cite{smale}. A detailed proof is offered by J. Milnor in \cite{h-cob} where the cancellation 
theorem of Morse occupies one chapter ({\it i.e.} Section 5) and is told to be {\it ``quite formidable''}.

Actually, Milnor's proof deals with the cancellation of two zeroes of a gradient (or {\it pseudo-gradient}) allowing the deformation to run among vector fields which {\it a priori} are not gradients. The effective deformation of the function remains hidden and the crossing of the 
codimension one stratum of functions with one cubic singularity is invisible. Our goal is to 
make this deformation visible; in particular the support of the deformation will be specified,
as it is in the notes by J. Cerf \& A. Gramain \cite{cerf} 
 (actually, we specify a smaller support 
than in \cite{cerf} where some {\it saturation} is not useful).
Moreover, we intend to show that the question of cancellation reduces to the one-dimensional
case where it is easy to solve. The same technique was already used in \cite{lauden}.

The statement we give here is slightly different  from the one in \cite{h-cob} which will be
derived as a
corollary. 
Moreover, we work with {\it descending pseudo-gradient}, called {\it pseudo-gradient} for 
short. According to K. Meyer \cite{meyer}, it is a vector field $X$ satisfying the two following conditions:
\begin{itemize}
\item[--] {Lyapunov inequality}: $X\cdot f<0$ apart from the critical points of $f$;
\item[--] Non-degeneracy condition: $X\cdot f$ has  a non-degenerate maximum
 at each critical point of $f$.
\end{itemize}
Such a pseudo-gradient $X$ may be built by partition of unity.
It follows  from the definition 
that if $p$ is a critical point of $f$, it is a hyperbolic zero of $X$ and hence, there are 
stable and unstable manifolds, respectively denoted by $W^s(p)$ and $W^u(p)$, which are  
formed by the points  $x\in M$ such that $X^t(x)$ tends to $p$ as $t$ tends to $+\infty$
or $-\infty$; here $X^t $ denotes the flow of $X$. The dimension of the unstable manifold $W^u(p)$
and the codimension of the stable manifold $W^s(p)$ both 
equal the index of the critical point $p$. It is easily  checked that the restriction of $f$ to $W^u(p)$
(resp. $W^s(p)$)
 has a non-degenerate maximum (resp. minimum).\\

\nd {\sc Theorem.} {\it Let us consider the Morse function $f:M^n\to \R$ 
equipped with a pseudo-gradient $X$. Let $p$ and $q$ be two critical points of $f$ satisfying the following conditions:

\nd {\rm i)} $W^u(p)$ and $W^s(q)$ intersect transversely
 and the intersection is made of one orbit $\ell$ of $X$ only.%
 
\nd  {\rm ii)} For  some $\ep>0$, each orbit of $X$ 
in $W^u(p)$ distinct from $\ell$ crosses the level
 set $f^{-1}\bigl(f(q)-\ep\bigr)$.

Then the pair $(p,q)$ is \emph{cancelable}. More precisely, if $U$ denotes an open neighborhood of
the closure of $W^u(p)\cap \{f\geq f(q)-\ep\}$, there is a path of 
smooth functions $\left(f_t\right)_{ t\in [0,1]},$ such that:

\nd {\rm 1)} $f_0=f$; 

\nd {\rm 2)} for every $t\in [0,1]$, $f_t$ coincides with $f$ on $M\smallsetminus U$;

\nd {\rm 3)}  the function $f_t\vert U$ is  Morse with two critical points  when $0\leq t<1/2$;
it has a cubic singularity when $t=1/2$ and it has no critical point when $1/2<t\leq 1$. \\
}

By computing the dimensions one has  the following formula for the Morse indices:
${\rm ind}(p)={\rm ind}(q)+1.$ 
There are other statements by varying assumption ii) which can be derived from the above theorem.

\section{Proof of the theorem}


 Say index$(p)=k+1={\rm index}(q)+1$. 
 We are given an open neighborhood $U$ of the closure 
of $W^u(p)\cap \{f(q)\leq f\leq f(p)\}$. After looking at $U\cap \{f=f(q)\}$, we may assume that the  frontier 
of $U$ in $M$ traces a cylinder $C$ tangent to $X$
in the domain $\{f(q)-\ep'\leq f\leq f(q)+\ep'\}$, where $\ep'>0$  satisfies assumption ii) of our theorem.
We then choose Morse charts $U(p)$ and $U(q)$ about $p$ and $q$ respectively,
 both contained in $U$.
They are equipped with coordinates $(\tilde y, z)\in\R^{k+1}\times \R^{n-k-1}$ 
(resp. $(y, \tilde z)\in \R^{k}\times \R^{n-k}$) such that
$f\vert U(p)= f(p)-\vert\tilde y\vert^2+\vert z\vert^2$ 
(resp. $f\vert U(q)= f(q)-\vert y\vert^2+\vert \tilde z\vert^2$). Moreover,
 it is possible to require that 
$(\tilde y, 0)$ are coordinates of $W^u(p)\cap U(p)$  (thus, $z$ must be  coordinates in the  
orthogonal complement to $T_pW^u(p)$ with respect to the Hessian of $f)$. Similarly, the coordinates 
of $U(q)$  may
 be chosen so that $(0,\tilde z)$ are  coordinates of $W^s (q)\cap U(q)$. 
 
 For $0< \ep < \ep'$  small enough there is a  Morse model $\M(q)\subset U(q)$ disjoint from $C$
whose top and bottom 
are  respectively at level $f(q)+\ep$ and $f(q)-\ep$;
and also,  there is a Morse model $\M(p)\subset U(p)$ with top and bottom at level
$f(p)+\ep$ and  $f(p)-\ep$.

 Denote by 
$\tilde \ell$ the trace of $\ell$ in the domain $\{f(q)+\ep\leq f\leq f(p)-\ep\}$. 
One end point $a$ of $\tilde\ell$
is $W^u(p)\cap\{f=f(p)-\ep\}\cap W^s(q)$  and the other end point $b$ is 
$W^u(p)\cap \{f=f(q)+\ep\}\cap W^s(q)$.  Consider the  diameter $\De_a\subset U(p)$ 
of the ball $\{\vert\tilde y\vert^2\leq\ep, z=0\}$
ending at $a$. Split the coordinates $(\tilde y, z)$ as $(t,y,z)$
 so that  $t$  is a coordinate of 
$\De_a$ and 
$$f\vert U(p)= f(p)-t^2 -\vert y\vert^2+\vert z\vert^2.$$  
 Similarly, 
  split $(y, \tilde z)= (y,t,z)$ in $U(q)$ so that $b\in \De_b:=\{y=0, z=0\}$ and 
  $$f\vert U(q)= f(q)-\vert y\vert^2
  +t^2 +\vert z\vert^2.$$ 

It is easy to find  
 a $C^\infty$ parametrized arc $\al : [0,1]\to M$, $u\in [0,1] \mapsto \al(u)\in M$,
 whose image is denoted by $A$,  and a sequence $0<u_0<u'_0<u'_1<u_1< 1$
with the following properties:
\begin{itemize}
\item   $p=\al(u_0)$, $a=\al(u'_0)$,
 $b=\al(u'_1)$, $q=\al(u_1)$ and $\al([u'_0,u'_1])=\tilde\ell$;
\item the image  $\al([0,u'_1])$ is contained in 
 $W^u(p)$ and the image  
$\al([u'_1,1])$ is contained in 
$W^s(q)$; 
\item   the function 
$h$ defined by $h:=f\circ\al$ has two critical points, one  maximum in $p$, one local minimum
in $q$ (both non-degenerate) and varies from $h(0)=f(q)-\ep$   to $h(1)= f(q)+\ep$.
\end{itemize}
One observes that there exists an increasing  smooth function $h_1$
without critical points,  which coincides with $h$
near the extremities of $[0,1]$, and satisfies 
$$(*)\quad\quad h_1\leq h\,. 
$$
This is the solution of the cancellation problem in the one-dimensional   case. \\

We now construct a $(k+1)$-dimensional sub-manifold $W\subset U$  containing $A$ such that 
$p$ and $q$ are the only critical points of 
$f\vert W$. This $W$ will contain the part of $W^u(p)$ above the level $f=f(q)+\ep$ that is,
$W^u(p)\cap\{f\geq f(q)+\ep\}\cong D^{k+1}$. Let $S$ be the boundary of this disc. Let $T$ denote the top
of $\M(q)$; we have $T\cong D^k\times S^{n-k-1}$ and the product structure of $T$ is determined 
by the Morse coordinates. The point $b$ belongs to the so-called {\it belt sphere} $\{0\}\times S^{n-k-1}
\subset W^s(q)$.
The radius of $D^k$ measured with the norm of the $y$-coordinates is a free choice for the Morse model
$\M(q)$.
If this radius is small enough, the transversality assumption implies that $S\cap T$ is isotopic 
in $T$ to the {\it standard meridian} $D^k\times \{b_0\}$ by an isotopy keeping the belt sphere fixed.
Therefore, by a first 
 level preserving isotopy of the coordinates of $\M(q)$ one makes $S\cap T$  be contained in $\{z=0\}$.
 By a second isotopy of this type
   one makes $X$ and $\nabla f$ coincide near the top of $\M(q)$; here $\nabla f$
stands for the descending gradient of $f$ in   $\M(q)$ associated to the Euclidean metric of the coordinates
$(t,y,z)$. These two isotopies  leave $W^s(q)$ invariant and are stationary near $q$.

Let 
$B\subset \M(q)$ be defined by $B:=\{z=0\}$; it is $C^\infty$ tangent to $W^u(p)$ along $B\cap S$.
 By partition of unity one finds
 a pseudo-gradient  $\xi$ such that:
\begin{itemize}
\item $\xi=X$ on $W^u(p)$ above the level $f=f(q)+\ep$ and on the cylinder $C$ defined previously;
\item $\xi=\nabla f$ in $\M(q)$.
\end{itemize}
Let $W^u_\xi (p)$ be the unstable manifold of $p$ for $\xi$. Each orbit of $\xi$ in $W^u_\xi (p)$
which does not meet $\tilde\ell$ reaches the level set $\{f=f(q)-\ep\}$ and remains  in $U$ since 
$C$ prevents it from  getting out. We define $W$
as the union of $B$ and $W^u_\xi (p)$ truncated by removing $\{f<f(q)-\ep\}$.\\

\nd {\sc Lemma 1.} {\it There exist a tubular neighborhood $N$ of $A$ in $U$
and coordinates $(u,y,z)\in [0,1]\times \R^k\times \R^{n-k-1}$ of  $N$, such that:
\begin{itemize}
\item[{\rm i)}] $A= \{y=0,z=0\}$;
\item[{\rm ii)}] $N_0:=N\cap W= \{z=0\}$;
 \item[{\rm iii)}] $f(u,y,z)= h(u)-\vert y\vert^2+\vert z\vert^2$ at every point of $N$.\\
\end{itemize}}

\nd {\sc Proof.} The first two items are easy to realize: take a provisory pair of small tubular neighborhoods 
$(N, N_0)$ of $A$ in $(M,W)$ whose fibres are contained in $\{t= const\}$ with respect to
 the coordinates of $\M(p)$ and $\M(q)$. Choose the coordinates $y$ in the fibre of $N_0$
 over $A$ and $z$ in the 
 fibre of $N$ over $N_0$  such that $(y,z)$ are the so-named coordinates in $\M(p)$ and $\M(q)$.
 For proving item iii), we need to invoke a general fact stated below.\\

 \nd {\sc Claim.} {\it Let $V$ be a smooth manifold and $V'\subset V$ be a sub-manifold. 
 Two germs of smooth functions $f$
and $g$
along  $V'$ 
whose restrictions  to $V'$  coincide and have no critical points
are isotopic  relative to $V'$ (meaning that they become equal after moving one of them by an isotopy of $V$
fixing $V'$); moreover, if $f=g$ near a compact set $K\subset V'$, the isotopy may be the identity
near $K$.} \\
 
For proving this claim the path method of Moser is available; it may  be applied to the path 
$t\in [0,1]\mapsto f_t:= (1-t)f+tg$, that is: there is  a time-dependent vector field $Z_t$ whose integration solves the problem of conjugating each $f_t$ to $f_0$ near $V'$ smoothly in $t\in[0,1]$. It is sufficient 
to find  $Z_t$ satisfying 
$df_t(x).Z_t(x)= f(x)-g(x)$ near $V'$ (see \cite{moser}); for that, apply the implicit function theorem
for finding local solutions which  are then glued by partition of unity. \bull

 We are going to apply this claim twice. First, apply it  to the data:  
 $V=W$,  
$V'=A$, $K= A\cap \bigl(\M(p)\cup \M(q)\bigr)$, the restriction $f\vert N_0$ and $g(u, y)=h(u)-\vert y\vert^2$.
 The isotopy yielded by the claim moves the coordinates $(u,y)$, without changing $u$ on $A$,
 so that the equality of iii) holds true in a neighborhood $N_1\subset N_0$ of $A$ in  $W$. 
 Then, apply the claim
 to the data: $V=M$, $V'=N_1$, $K= N_1\cap  \bigl(\M(p)\cup \M(q)\bigr)$,
 $f$ and $g(u,y,z)= h(u)-\vert y\vert^2+\vert z\vert^2$. 
The new isotopy finishes to put the coordinates in a position which makes the equality of  item iii)
hold true near $A$ in $M$ yielding the desired $N$.\bull\\

In $W\cap N$,
the vector field $Y:=y\partial_y$  is tangent to $W$ and  it is   a pseudo-gradient  for 
$f\vert \bigl(W\cap (N\smallsetminus A)\bigr)$. Still denoted 
$Y$, it extends to the whole $W$
as a pseudo-gradient  for $f\vert (W\smallsetminus A)$ (by partition of unity). Each  orbit  of 
$Y$ converge to a point of $A$ in the past. This defines a  $k$-disk fibration $W\to A$
which is pinched at $\al(0)\in A$. Indeed, this fibration 
is clear for $W\cap N \to A$ and  the orbits of $Y$
yield an ambient isotopy  from $W$ to $W\cap N$. Denote $D_u$ the fibre of $W$ 
over $\al(u)$.

Let $\widetilde W$ be a tubular neighborhood of $W$ in $\{f\geq f(q)-\ep\}$. The fibration $W\to A$
extends to $\widetilde W$ as a $(n-1)$-disk fibration:
$$\widetilde D_u \hookrightarrow \widetilde W\to A\,,
$$
pinched at $\al(0)$ and coinciding with the projection $(u,y,z)\to \al(u)$ in $N$. The function
 $ f_u:= f\vert \widetilde D_u$ 
 is a Morse function
with $\al(u)$ as unique critical point and its  index is  $k$. It has a pseudo-gradient
$\widetilde Y_u$, tangent to $\widetilde D_u$, which coincide with $Y$ on $W$ 
and with $y\partial_y-z\partial_z$ on $N$.\\

\nd{\sc Lemma 2} (Decreasing of a critical value.) {\it  Let $g:V\to\R$ be a Morse function with a pseudo-gradient $Z$. Let $p$
be a critical point of index $k$ and let $a<g(p)$. Assume that the unstable manifold
 $W^u(p)$ contains a compact $k$-disk  $D$
whose  boundary lies  in the level set $\{g=a\}$.  Let $U$ be a neighborhood of $D$ in $\{g\geq a\}.
$Then, for every $\ep>0$,
there exists a one parameter family of Morse functions $\left(g_t\right)_{t\in[0,1]}$,
such that:\begin{itemize}
\item $g_0=g$, $a<g_1(p)<a+\ep$,
\item $Z$ is a pseudo-gradient of $g_t$ for every $t\in[0,1]$,
\item and $g_t=g$ out of $U$.
 \end{itemize}
 The same statement holds true with parameters.\\
}

\nd {\sc Proof.}  The point $p$ is a hyperbolic zero of the vector field $Z$. Then, the orbits
of $Z$ close to $D$ but not tangent to $D$ are crossing the level set $\{g=g(p)+\eta\}$, for a small
$\eta>0$ and, of course,  the level set $\{g=a\}$. Consider the foliation $\mathcal F$
defined by $g=const$
on $U\smallsetminus D$.  Find an $n$-dimensional compact domain
 $K$ in $U\smallsetminus D$ made of pieces of $Z$-orbits from $\{g=g(p)+\eta\}$ to $\{g=a\}$ 
 surrounding $D$. It is easy to replace $\mathcal F\vert K$
 by a foliation whose leaves are still transverse to $Z$ in order to have the level sets of the wanted function $g_1$. We refer to
\cite{lauden} for more details (see also \cite{h-cob} Section 4 or
 \cite{{smale_gr}} Section 2).\bull\\

\nd{\sc End of the proof of the theorem}. The critical point of $f_u$ is $p_u:=\al(u)=(u,0,0)$ and its value is $h(u)$.
The vector field $\widetilde Y_u$ is a pseudo-gradient for $f_u$
and the disk $D_u$ fulfills the requirement
of Lemma 2, which we apply smoothly with respect to  the parameter $u$.
 So, it is possible to decrease the critical value of $f_u$,
 smoothly in $u\in[0,1]$, from $h(u)$ to $h_1(u)$, still keeping $\widetilde Y_u$ as a pseudo-gradient.
Therefore, the deformation of $f$ over  $A$ 
from $h$ to $h_1$ (compare $(*)$),
extends to $\widetilde W$ without creating new critical points
 away from $A$.
Moreover, the deformation is supported in $\widetilde W$, hence in $U$, as wanted.\bull 

\section{Applications}
In this section, there are given two  more classical statements, starting with the one from 
Milnor's book (\cite{h-cob}, Section 5).\\

\nd {\sc Corollary 1.} {\it Let $(W,L_0,L_1)$ be a compact cobordism 
and $f: (W,L_0,L_1)\to
 ([0,1],0,1)$ be a Morse function with two critical points $p$ and $q$. Assume there is a pseudo-gradient such that $W^u(p)$ and $W^s(q)$ intersect in one orbit and transversely.
 Then the cobordism is a product: $W\cong M_0\times[0,1]$.}\\
 
 
 \nd {\sc Proof.} Certainly $0<f(q)<f(p)<1$. 
 By compactness,
 every orbit of the pseudo-gradient in $W^u(p)$ distinct from the connecting orbit
 reaches $L_0$. Hence, the assumptions of our theorem are fulfilled. After canceling the critical points, $W$ is a product by gradient lines.\bull\\
 
 \nd {\sc Corollary 2.} {\it Let $M$ be a closed manifold, $f:M\to\R$ be a Morse function, $(p,q)$
 be a pair of critical points whose respective indices are $k+1$ and $k$.
  Let $X$ be a pseudo-gradient of $f$. Assume that there are  compact disks  
  $D(p)\subset W^u(p)$ and $D(q)\subset W^s(q)$ with the following properties:
  \begin{itemize}
  \item $D(p)$ and  $D(q)$ are neighborhoods   of $p$ and $q$  respectively in $W^u(p)$ and 
  $W^s(q)$);
  \item their boundaries lie in a  regular 
  level set $\{f=a\}$, $f(q)<a<f(p)$ in which they intersect in one point
  only and transversely.
  \end{itemize}
  Then, the pair $(p,q)$ is cancelable. }\\
  
  \nd {\sc Proof.} Since $a$ is a regular value, the disk $D(q)$ can be extended (keeping its name)
  so that its boundary lies in $\{f=a+\ep\}$, for some small $\ep>0$. Lemma 2 can be applied
  for decreasing the critical value $f(p)$ so that $a<f(p)<a+\ep$.
 Similarly, by considering $-f$,
   Lemma 2 can be applied 
  for increasing the critical value $f(q)$ so that $a<f(q)<f(p)<a+\ep$. In this situation, as in Corollary 1, the assumptions of our theorem are fulfilled. The cancellation follows.\bull\\


%% file: cancel_hal.bbl
\vskip .5cm
\begin{thebibliography}{99}

\bibitem{cerf} J. Cerf, A. Gramain, {\it Le th\'eor\`eme du
 $h$-cobordisme (Smale)}, cours Orsay 1966, Secr. math. \'Ec. Normale Sup., Paris, 1968  (www.maths.ed.ac.uk/~aar/surgery/cerf-gramain.pdf).


\bibitem{hueb} W. Huebsch, M. Morse, {\it The bowl theorem and a model non-degenerate function}, Proc. Nat. Acad. Sc. U.S.A., vol. 51 (1964), 49-51.
\bibitem{lauden} F. Laudenbach, {\it A proof of Reidemeister-Singer's theorem by Cerf's methods},
arXiv(math-GT): 1202.1130.

\bibitem{meyer} K. Meyer, {\it Energy functions for Morse-Smale systems}, Amer. J. Math. 90 (1968), 1031-1040.
\bibitem{h-cob} J. Milnor, {\it  Lectures on the h-cobordism theorem,} 
Princeton Univ. Press, 1965.


\bibitem{morse1} M. Morse, {\it The existence of polar non-degenerate functions
on differentiable manifolds}, Annals of Math. 71 (1960), 352-383.

\bibitem{morse2} M. Morse, {\it Bowls of a non-degenerate function on a compact differentiable 
manifold}, 81-103 in: Differential and Combinatorial Topology (A Symposium in Honor of Marston Morse), Princeton Univ. Press, 1965.

\bibitem{moser} J. Moser, {\it On the volume elements on a manifold,} Trans. Amer. Math. Soc.
120 (1965), 286-294.

\bibitem{smale_gr} S. Smale, {\it On gradient dynamical systems,} Annals of Math.
74 (1961), 199-206.
\bibitem{smale} S. Smale, {\it Generalized Poincar\'e's conjecture in dimensions greater than four}, Annals of math. 74 (2) (1961), 391-406.
\bibitem{thom49}    R. Thom, {\it Sur une partition en cellules associ\'ee \`a une fonction sur une vari\'et\'e}, Comptes Rendus Acad. Sc. Paris, t. 228, (1949), 973-975.


 \end{thebibliography}
